\documentclass[11pt,oneside]{amsart}
\usepackage{amsmath,amscd}
\usepackage{ifthen}
\usepackage{amsfonts}
\usepackage{amssymb}
\usepackage{epsf}
\usepackage{amsopn}

\newcommand{\showcomments}{yes}

\newcommand{\comment}[1]
{\ifthenelse{\equal{\showcomments}{yes}}
{\footnotemark\marginpar{\sffamily{\tiny
\addtocounter{footnote}{-1}\footnotemark#1

}\normalfont}}{}}

\newtheorem{thm}{Theorem}[section]
\newtheorem{lem}[thm]{Lemma}
\newtheorem{cor}[thm]{Corollary}

\newtheorem*{torsion rf theorem}{Theorem~\ref{thm:torsion rf}}

\theoremstyle{definition}
\newtheorem{defn}[thm]{Definition}
\newtheorem{rem}[thm]{Remark}
\newtheorem{exmp}[thm]{Example}

\newtheorem{conv}[thm]{Convention}

\newcommand{\refe}[1]{\eqref{e:#1}}

\begin{document}

\author[Inna Bumagin]{Inna Bumagin}
\address{School of Mathematics and Statistics, Carleton University,  1125 Colonel By Drive, Ottawa, ON, Canada, K1S 5B6}

\email{bumagin@math.carleton.ca}

\author[O. Kharlampovich]{Olga Kharlampovich}
\address{Department of Mathematics and Statistics, McGill University, 805 Sherbrooke Street West, Montreal, QC, Canada, H3A2K6}

\email{olga@math.mcgill.ca}
\thanks{The first and the second author were supported
 by NSERC Grant}

\author[A. Miasnikov]{Alexei Miasnikov}
\address{Department of Mathematics and Statistics, McGill University, 805 Sherbrooke Street West, Montreal, QC, Canada, H3A2K6}
\email{alexeim@att.net}
\thanks{The third author was supported by NSERC Grant and by NSF
GrantDMS-9970618}

\title[Isomorphism problem]{Isomorphism problem for finitely generated fully residually free groups}



\begin{abstract} We prove that the isomorphism problem for finitely generated fully
residually free groups (or $\mathcal{F}$-groups for short) is
decidable. We also show that each $\mathcal{F}$-group $G$ has a
decomposition that is invariant under automorphisms of $G$, and
obtain a structure theorem for the group of outer automorphisms
$Out(G)$.
\end{abstract}

\maketitle
\tableofcontents
\section{Introduction}
The isomorphism problem - find an algorithm that for any two
finite presentations determines, whether or not the groups defined
by these presentations are isomorphic - is the hardest of the
three algorithmic problems in group theory formulated by Max Dehn
at the beginning of the 20th century. It is easy to see that
solvability of the isomorphism problem in the class of finitely
presented groups implies solvability of the word problem (find an
algorithm to determine, whether or not a given product of
generators of a group represents the trivial element of the
group). The isomorphism problem is unsolvable in the entire class
of finitely presented groups, because there exist finitely
presented groups with unsolvable word problem; this latter
assertion is the fundamental result of Novikov and Boone. One can
still try to solve the isomorphism problem restricted to a certain
class $\mathcal{C}$ of finitely presented groups: \textit{find an
algorithm that for any two finite presentations of groups from the
class $\mathcal{C}$ determines, whether or not the groups defined
by these presentations are isomorphic}. There are only few classes
of groups for which the isomorphism problem is known to be
solvable. This is a classical result that the isomorphism problem
is solvable for finitely generated Abelian groups. Solvability of
the isomorphism problem for finitely generated free groups has
been known since 1950ties due to the work of Nielsen. Among the
most significant results in this area is Segal's solution to the
isomorphism problem for polycyclic-by-finite groups~\cite{Segal}.
One should also mention the positive solution to the isomorphism
problem for finitely generated nilpotent groups, which is an
earlier result obtained by Segal and Grunewald~\cite{SeGr}.
Another profound result was obtained by Sela~\cite{Se2} who proved
that the isomorphism problem is solvable for torsion-free word
hyperbolic groups which do not split over a cyclic subgroup. One
of the most important ingredients of Sela's solution to the
isomorphism problem is the decidability of equations over free
groups proved by Makanin~\cite{Ma} and Razborov~\cite{Ra}, and
extended by Rips and Sela~\cite{RS1} to torsion-free word
hyperbolic groups.

We consider the class of \emph{finitely generated fully residually
free groups} ($\mathcal{F}$-groups for short) defined as follows.
\begin{defn}\label{defn:fullyresfree}\cite{BB} A group $G$ is called \emph{fully
residually free} if for any finite number of non-trivial elements
$g_1,\dots,g_n$ in $G$ there exists a homomorphism $G\rightarrow
F$ from $G$ to a free group $F$ that maps $g_1,\dots,g_n$ to
non-trivial elements of $F$.
\end{defn}
The first examples of non-free fully residually free groups are
due to Lyndon~\cite{Lfrf}, where he introduced free Lyndon's
$\mathbb{Z}[t]$-groups and proved that they are fully residually
free. In the same year 1960, in a very influential
paper~\cite{Leq} he used these groups to describe completely the
solution sets of one-variable equations over free groups.

A finitely generated fully residually free group $G$ is word
hyperbolic, if any maximal Abelian subgroup of $G$ is
cyclic~\cite{KM2}. However, in this latter case $G$ has one of the
following decompositions: a non-trivial free decomposition, or a
non-trivial JSJ decomposition, or $G$ is the fundamental group of a
closed surface and has a non-trivial cyclic splitting. Therefore,
the case of a word hyperbolic fully residually free group is not
covered by Sela's solution to the isomorphism problem. Our main
result is the following theorem.
\setcounter{section}{4}
\setcounter{thm}{12}
\begin{thm} Let $G\cong\langle
\mathcal{S}_G\mid\mathcal{R}_G\rangle$ and $H\cong\langle
\mathcal{S}_H\mid\mathcal{R}_H\rangle$ be finite presentations of
fully residually free groups. There exists an algorithm that
determines whether or not $G$ and $H$ are isomorphic. If the
groups are isomorphic, then the algorithm finds an isomorphism
$G\rightarrow H$.
\end{thm}
\setcounter{section}{1}\setcounter{thm}{1}
The most important ingredients of our proof are computability of a
JSJ decomposition of an $\mathcal{F}$-group, and solvability and
the structure of the solution sets of equations over
$\mathcal{F}$-groups, obtained by the second and the third
authors~\cite{KM3},~\cite{JSJ} (see also
Theorem~\ref{thm:effectiveJSJ} and Section~\ref{s:algorithm} in
the present paper). To deduce solvability of the isomorphism
problem, we prove that a one-ended $\mathcal{F}$-group $G$ has a
canonical Abelian JSJ decomposition that is invariant under
automorphisms of $G$. Moreover, using results obtained
in~\cite{JSJ}, we deduce that the canonical decomposition can be
constructed effectively. More precisely, in
Theorem~\ref{thm:universalJSJ} we define an Abelian JSJ
decomposition $\Gamma(V,E)$ of $G$ that has the following
property.
\begin{thm}\label{thm:canonicalJSJ} Let $G$ be a one-ended
$\mathcal{F}$-group, and let $\Gamma(V,E)$ be an Abelian JSJ
decomposition of $G$ that satisfies the conditions of
Theorem~\ref{thm:universalJSJ}. If a graph of groups $\Delta(U,P)$
is another Abelian JSJ decomposition of $G$ that satisfies the
conditions of Theorem~\ref{thm:universalJSJ} also, then $\Delta$
can be obtained from $\Gamma$ by conjugation and modifying
boundary monomorphisms.
\end{thm}
Theorem~\ref{thm:canonicalJSJ} follows from
Theorem~\ref{thm:down_to_vertices}. Hyperbolic groups have
canonical JSJ decompositions over virtually cyclic subgroups as
was shown by Bowditch~\cite{Bow}, this result was first proved by
Sela~\cite{Se2} for torsion-free hyperbolic groups. Another class
of groups that possess canonical JSJ decompositions was introduced
by Forester~\cite{For} (Guirardel~\cite{Guir} gave an alternate
proof of this latter result). Not all finitely presented groups
have canonical JSJ decompositions, as shown by
Forester~\cite{Foruniq}. Using Theorem~\ref{thm:canonicalJSJ}, we
obtain the following structure theorem for $Out(G)$ (cf.
Theorem~\ref{thm:auto group}).
\begin{thm}\label{thm:auto group intro} Let $G$ be a one-ended
$\mathcal{F}$-group. $Out(G)$ is virtually a direct product of a
finitely generated free Abelian group, subgroups of
$GL_n(\mathbb{Z})$, and the quotient of a direct product of
mapping class groups of surfaces with boundary by a central
subgroup isomorphic to a finitely generated free Abelian group.
\end{thm}
Recall that similar results for torsion-free hyperbolic groups
were obtained by Sela~\cite{Sela} and for a more general class of
groups by Levitt~\cite[Theorem~1.2]{Levitt}.

The first author wishes to thank Ilya Rips, Zlil Sela and Daniel
Wise for numerous useful conversations preceding the work on the
present paper.

\section{Graphs of groups and splittings}\label{s:def}
\begin{defn}\label{defn:graph} A {\em directed graph} $(V,E,\mathcal{O})$
consists of a set of vertices $V$, a set of edges $E$ and an
orientation $\mathcal{O}$ determined by two functions $i\colon E
\rightarrow V$ and $\tau \colon E \rightarrow V$. For an edge
$e\in E$ the vertex $i(e)$ is the \emph{initial vertex} of $e$,
and $\tau (e)$ is the {\it terminal vertex} of $e$. We call $i(e)$
and $\tau(e)$ the \emph{endpoints of e}.
\end{defn}
\begin{defn}\label{defn:graphgps} A graph of groups $\Gamma(V,E,\mathcal{O})$
is a directed graph $(V,E,\mathcal{O})$ where to each vertex $v\in
V$ (or to each edge $e\in E$) we assign a group called $G_v$ (or
$G_e$) so that for each edge $e\in E$ there are monomorphisms
\[
\alpha\colon G_e\rightarrow G_{i(e)}\quad \text{and}\quad
\omega\colon G_e\rightarrow G_{\tau(e)}
\]
called the \emph{boundary monomorphisms from the edge group} $G_e$
to the \emph{vertex groups} $G_{i(e)}$ and $G_{\tau(e)}$. We refer
to $G_v$ and $G_e$ the \emph{stabilizer} of $v$ and $e$,
respectively.
\end{defn}

\begin{defn}\label{defn:splitting} By a \emph{splitting of} $G$ we mean
a \emph{triple} $\Sigma=(\Gamma(V,E,\mathcal{O}),T,\varphi)$ where
$\Gamma(V,E,\mathcal{O})$ is a graph of groups, $T$ is a
\emph{maximal subtree of the graph} $(V,E)$ and
$\varphi\colon\pi_1(\Gamma(V,E,\mathcal{O});T)\rightarrow G$ is an
\emph{isomorphism}.

We recall that the fundamental group of a graph of groups
$\pi_1(\Gamma(V,E,\mathcal{O});T)$ with respect to a maximal
subtree $T$ is given by
\begin{align*}
\langle G_v (v\in V), t_e (e\in E_0)\mid &\forall e\in E_0
(t_et_{\bar{e}}=1 , \alpha(g)t_e=t_e\omega(g), \forall g\in G_e),\\
&\forall e\in T (\alpha(g)=\omega(g),\forall g\in G_e)\rangle
\end{align*}
where $E_0=\{e\in E\mid e\notin T\}$ denotes the set of edges that
do not belong to the maximal tree.
\end{defn}
Let $G$ be a group and let $\mathcal{G}$ be a set of splittings of
$G$ into a graph of groups.  One introduces an equivalence
relation on $\mathcal{G}$ generated by the following operations
(we refer the reader to~\cite{RS2} and to~\cite[Section~2.4]{JSJ}
for more details):
\begin{enumerate}
  \item \emph{Conjugation} is a usual conjugation;
  \item \label{e:modify} \emph{Modifying boundary monomorphisms by conjugation} is
  defined as follows. Let $G=\langle A,t\mid
t\alpha(c)t^{-1}=\omega(c)\forall c\in
  C\rangle$. For an arbitrary element $a\in A$ one defines
$\alpha'\colon C\rightarrow
  A$ by $\alpha'(c)=a^{-1}\alpha(c)a$, and replaces $\alpha$ by
  $\alpha'$. One replaces also the isomorphism $\varphi$ by the
  isomorphism $\varphi_a$ defined by $\varphi_a(t)=\varphi(ta)$
  and $\varphi_a(g)=\varphi(g)$ for all $g\neq t$. If $G=A\ast_C B$,
  then one replaces the monomorphism
  $\alpha\colon C\rightarrow A$ by $\alpha'\colon C\rightarrow A$
  defined as above and $\varphi$ by the isomorphism $\varphi_a$
  defined by $\varphi_a(g)=\varphi(g)$ for $g\in A$
  and $\varphi_a(g)=\varphi(a^{-1}ga)$ for all $g\in B$. For a
  general graph of groups, let $e$ be the edge stabilized by $C$;
  one collapses all edges but $e$ and defines $\alpha'$ and
  $\varphi_a$ as above, with the only restriction that $a\in
  G_{i(e)}$.
  \item \emph{Sliding} corresponds to the relation $$(A\ast_{C_1}
B)\ast_{C_2}
  D\cong(A\ast_{C_1}D)\ast_{C_2}B$$ in the case when $C_1\subseteq
C_2$.
 \item\label{e:op4} By a \emph{refinement of $\Delta\in\mathcal{G}$ at a vertex}
$v\in\Delta$ we mean replacing $v$ by a non-degenerate graph of
groups $\gamma(V_{\gamma},E_{\gamma})$ which is compatible with
$\Delta$ and has the fundamental group $G_v$ (where $G_v$ is the
stabilizer of $v$ in $\Delta$). A vertex $v$ is {\it flexible} if
there exists a refinement of $\Delta$ at $v$; otherwise, $v$ is
{\it rigid}.
\end{enumerate}
In what follows, by a splitting of $G$ we mean a graph of groups
$\Gamma(V,E)$; when there is no ambiguity, we identify the groups
assigned to edges $G_e$ with their images $\alpha(G_e)\subseteq
G_{i(e)}$ and the groups assigned to vertices with their images in
$G$ under the isomorphism $\varphi$. Usually, we do not specify a
maximal tree and an orientation in the graph $(V,E)$. Observe that
conjugation corresponds to an inner automorphism of $G$, whereas
operations~\refe{modify}-~\refe{op4} change the presentation of $G$
as a graph of groups and usually do not lead to an automorphism of
$G$. However, there is an exception. If we have an operation of
type~\refe{modify} so that $a$ is in the centralizer
$C_A(\alpha(C))$ of $\alpha(C)$ in $A$, then $\alpha'(C)=\alpha(C)$
which means that we actually do not modify the graph of groups.
Then, in the above notation, the composition of the isomorphisms
$\varphi_a\circ\varphi^{-1}$ is well-defined and results in an
automorphism of $G$ called a \emph{generalized Dehn twist}. More
precisely, we have the following definition.
\begin{defn}\label{defn:Dehn twist} Let $\Gamma(V,E)$ be an Abelian
splitting $G$ of a group $G$, and let $e\in E$ be an edge with the
endpoints $i(e)=v$ and $\tau(e)=u$ and the stabilizer $G_e$. By a
\emph{generalized Dehn twist along the edge} $e$ we mean an
automorphism $\beta_a\colon G\rightarrow G$ with $a\in
C_{G_v}(G_e)$, defined as follows.

If $e$ is a separating edge, let $\Delta_v$ (or $\Delta_u$) denote
the connected component of $(V,E)\setminus\{e\}$ that contains $v$
(or $u$). Then $\beta_a(g)=g$ for $g\in G_w$ with $w\in\Delta_v$
and $\beta_a(g)=aga^{-1}$ for $g\in G_w$ with $w\in\Delta_u$.

If $e$ is a non-separating edge, then one can choose a maximal
tree $T$ in $(V,E)$ so that $e\notin T$. Let $t$ be the stable
letter that corresponds to $e$. We set $\beta_a(t)=at$ and
$\beta_a(g)=g$ for all $g\neq t$.
\end{defn}
In particular, if the edge group $C$ is cyclic and
$\alpha(C)=C_A(\alpha(C))$, then our definition coincides with the
definition of a \emph{Dehn twist} (see~\cite{RS2}).
\begin{defn}\label{defn:kinds of split} A splitting is
\emph{elementary} if the graph $\Gamma(V,E)$ is either an edge of
groups or a loop of groups so that either $G\cong A*_{C}B$, or
$G\cong A*_{C}$. A splitting is called \emph{Abelian} if the edge
groups are all Abelian.
\end{defn}
\subsection{$G$-tree}
By a \emph{tree} we mean a simplicial tree i.e., a graph with no
circuits. One assigns unit length to each edge of a tree, to make
a tree into a geodesic metric space.
\begin{defn}\label{defn:Gtrees} A tree equipped with an action of a group $G$ is called
a $G$-\emph{tree}. An action $G\times X\rightarrow X$ is
\emph{Abelian}, if edge stabilizers in $X$ are all Abelian
subgroups of $G$. A $G$-tree $X$ is \emph{minimal} if it contains
no $G$-invariant proper subtrees. Two vertices (or edges) $x_1$
and $x_2$ in $X$ are $G$-\emph{equivalent}, if they belong to a
$G$-orbit.

By the \emph{fixed set of} $g\in G$ we mean $Fix(g)=\{x\in X\mid
g.x=x\}$. A $G$-tree is \emph{k-acylindrical}, if
$diam(Fix(g))\leq k$ for all $g\in G$.
\end{defn}
\begin{conv}\label{conv:abelian splittings} In what follows, we
consider Abelian splittings and Abelian actions, only.
\end{conv}
The central result of the Bass-Serre theory~\cite{STr},\cite{Bass}
tells that to each splitting $\Sigma=(\Gamma(V,E),T,\varphi)$ of a
group $G$ one can associate a minimal $G$-tree, which is the
covering space of the graph of groups $\Gamma(V,E)$, and vice
versa, $G$ inherits a splitting from its action on a minimal
$G$-tree with no inversions.
\subsection{Extended fundamental domain and natural lift}
\begin{defn}\label{defn:extended fund domain}
An \emph{extended fundamental domain} $D$ is a finite subtree of
$X$ so that the $G$-orbit of $D$ is the whole tree $X$, and
different edges of $D$ belong to different $G$-orbits.
\end{defn}
\begin{lem}\label{lem:ext domain} Vertices $v$ and $u\neq v$ of an
extended domain $D$ are $G$-equivalent if and only if either
$v=t.u$ or $u=t.v$, where $t$ is a stable letter in the
presentation of $G$ determined by $\Delta$.
\end{lem}
Proof of this lemma is straightforward and we omit it.
\begin{defn}\label{defn:reduced} A graph of groups $\Delta$ is
\emph{reduced}, if for each vertex $v$ of valency one or two,
$G_v$ properly contains the groups of adjacent edges. We say that
$\Delta$ is \emph{semi-reduced}, if for each edge $e\in E$ with
the endpoints $v$ and $u$, the equality $G_e=G_v$ implies that
$v\neq u$, $val(v)\geq 2$ and $G_e\varsubsetneq G_u$. We say that
a $G$-tree $X$ is \emph{(semi-)reduced}, if the corresponding
graph of groups $\Delta=G\backslash X$ is (semi-)reduced.
\end{defn}
\begin{defn}\label{defn:natural lift} Let $X$ be
$2$-acylindrical and semi-reduced. A \emph{natural lift $\lambda$
of $\Delta$ to }$X$ is defined as follows. The image of a vertex
$v\in \Delta$ with the stabilizer $G_v$ is the vertex
$\lambda(v)\in X$ with $Stab(\lambda(v))=G_v$. Let $e$ be an edge
with the endpoints $i(e)=v$ and $\tau(e)=u$. If $e\in T$, then
$\lambda(e)$ is the edge of $X$ joining $\lambda(v)$ and
$\lambda(u)$, and if $e\notin T$, then $\lambda(e)$ is the edge of
$X$ joining $\lambda(v)$ and $t_e.\lambda(u)$ where $t_e$ is the
stable letter corresponding to $e$.
\end{defn}
\begin{lem}\label{lem:lift basics}
\begin{enumerate}
  \item \label{e:lb1} The natural lift of $\Delta$ to $X$ is
well-defined.
  \item \label{e:lb2} The natural lift of $\Delta$ to $X$ is an
  extended domain.
\end{enumerate}
\end{lem}
\begin{proof} Assume that there are two vertices $x_1$ and $x_2$ in
$X$ with $Stab(x_1)=Stab(x_2)=G_v$. The path $p$ joining $x_1$ and
$x_2$ in $X$ is stabilized by $G_v$. Since $X$ is
$2$-acylindrical, the length of $p$ is either $1$ or $2$. If $p$
is an edge, we get a contradiction as $X$ is semi-reduced. Let the
length of $p$ equal $2$, and let $v=\pi(x_1)$ and $u=\pi(x_2)$ be
the natural projections of $x_1$ and $x_2$ to $\Delta$. Assume
that $val(v)>1$. The stabilizer of an edge $f\notin\pi(p)$
incident on $v$ is a non-trivial subgroup $B$ of $G_v$. The edge
$f\in\Delta$ lifts to an edge $q_f\in X$ so that $q_f\notin p$
with $Stab(q_f)=B$, so that the subgroup $B$ fixes $3$ edges in
$X$, a contradiction. Therefore, $val(v)=val(u)=1$ while $X$ is
semi-reduced, a contradiction. Thus, the image of each vertex in
$\Delta$ under a natural lift is defined uniquely. Since the
images of edges are determined uniquely by the images of their
endpoints, the assertion~\refe{lb1} follows. Furthermore, the
definition of the Bass-Serre tree $X$ as a covering space of
$\Delta$ implies the assertion~\refe{lb2}. Indeed, the $G$-orbit
of the natural lift of $\Delta$ is the whole $X$. Moreover, the
edges of $\Delta$ are representatives of different $G$-orbits of
edges in $X$, hence their lifts to $X$ are not $G$-equivalent.
%
\end{proof}

\subsection{Morphisms of graphs}
\begin{defn}\label{defn:iso of graphs} Let $(V,E)$ and $(U,B)$
be two graphs. A map $\chi\colon(V,E)\rightarrow(U,B)$ is
\emph{simplicial}, if $\chi$ maps each vertex $v\in V$ to a vertex
$u\in U$ and each edge $e\in E$ to a (possibly, empty) path in
$(U,B)$ so that the incidence relations are preserved. A
simplicial map $\chi\colon(V,E)\rightarrow(U,B)$ is an
\emph{isomorphism of graphs} if $\chi$ maps each edge $e\in E$ to
an edge $b\in B$ and is bijective on both the set of vertices and
the set of edges.
\end{defn}
\begin{rem} \label{rem:effective iso} It follows immediately from
the definition that for finite graphs $(V,E)$ and $(U,B)$ one can
find effectively the (possibly, empty) set of all isomorphisms
$\chi\colon(V,E)\rightarrow(U,B)$.
\end{rem}
\begin{defn}\label{defn:monomorphic image of JSJ} Let $\psi\colon G\rightarrow
H$ be an isomorphism of groups, and let $\mathcal{G}$ (or
$\mathcal{H}$) be the set of all splittings of $G$ (or $H$) into a
graph of groups. With the isomorphism $\psi$ we associate a map
$\psi_*\colon\mathcal{G}\rightarrow\mathcal{H}$, where the
\emph{image} $\psi_*(\Gamma)$ \emph{of}
$\Gamma(V,E)\in\mathcal{G}$ is the graph of groups
$\Delta(U,B)\in\mathcal{H}$ defined as follows:
\begin{enumerate}
  \item The underlying graphs $(V,E)$ and $(U,B)$ are isomorphic, and we
  identify each vertex and each edge of $(V,E)$ with its image in $(U,B)$ under an isomorphism.
  \item The group assigned to a vertex or to an edge in $\Delta$ is
  the $\psi$-image of the group assigned to that vertex or edge in $\Gamma$.
  \item Let $e$ be an edge with the endpoints $v=i(e)$ and $u=\tau(e)$.
  The boundary monomorphisms $\alpha_{\psi}\colon G_e\rightarrow G_v$ and
  $\omega_{\psi}\colon G_e\rightarrow G_u$ in $\Delta$ are defined by
$\alpha_{\psi}(\psi(b))=\psi(\alpha(b))$ and
$\omega_{\psi}(\psi(b))=\psi(\omega(b))$ for all $b\in G_e$.
\end{enumerate}
\end{defn}
\subsection{Universal decomposition of a group}\label{s:universal}
\begin{defn}\cite{BW}\label{defn:universal decomposition}
By a {\it universal decomposition} of $G$ we mean a decomposition
of $G$ into a graph of groups $\Gamma=\Gamma(V,E)$ that has the
following property. Given a minimal $G$-tree $T$, one can find
refinements at flexible vertices of $\Gamma$ and obtain a
decomposition $\Gamma_r$ of $G$ so that there exists a
$G$-equivariant simplicial map from the Bass-Serre tree
$\tilde{\Gamma}_r$ onto $T$.
\end{defn}
\begin{exmp} \label{ex:trivial decomp} Obviously, every group
$G$ has a \emph{trivial universal decomposition} that consists of
a unique flexible vertex stabilized by $G$. It can be readily seen
that if $G$ is a free (Abelian or non-Abelian) group or a closed
surface group, then in fact, the only universal decomposition of
$G$ is the trivial decomposition. More precisely, $G$ is
indecomposable in the sense of Definition~\ref{defn:indecomposable
group} below.
\end{exmp}
In what follows, we will be interested in an Abelian universal
decomposition of a group $G$ with maximal number of vertices. For
instance, the Grushko free decomposition is a maximal universal
decomposition in the class of all free decompositions of $G$. For
a freely indecomposable group, a JSJ decomposition has the
universal property (see Section~\ref{s:properties} for more
details).
\begin{defn}\label{defn:nondegen} We say that a graph of groups
$\Delta$ is \emph{non-degenerate}, if $\Delta$ is semi-reduced and
the set of edges of $\Delta$ is not empty.
\end{defn}
\begin{defn}\label{defn:indecomposable group} A group $G$ is
\emph{decomposable} if $G$ has a non-degenerate universal
decomposition. Otherwise, $G$ is \emph{indecomposable}. In
particular, if $G$ is an indecomposable group which is not a free
non-Abelian group, then $G$ is freely indecomposable.
\end{defn}

\section{Properties of fully residually free groups}\label{s:properties}
As before, we denote by $\mathcal{F}$ the \emph{class of finitely
generated fully residually free groups} (also called \emph{limit
groups} by Sela~\cite{Sd}), and say that $G$ is an
$\mathcal{F}$-group if $G$ belongs to the class $\mathcal{F}$. In
Theorem~\ref{thm:properties} below we mention only those
properties of $\mathcal{F}$-groups which we use in our proof.
\begin{thm}\label{thm:properties} Let $G$ be an
$\mathcal{F}$-group. Then $G$ possesses the following properties.
\begin{enumerate}
  \item \label{e:pr0} $G$ is torsion-free;
  \item \label{e:pr1} Each subgroup of $G$ is an
$\mathcal{F}$-group;
    \item \label{e:pr3} $G$ has the CSA property. Namely, each maximal Abelian
  subgroup of $G$ is malnormal, so that if $M$ is a maximal Abelian
  subgroup of $G$ then $M\cap gMg^{-1}\neq\{1\}$ for $g\in G$
  implies that $g\in M$;
  \item \label{e:pr2} Each Abelian subgroup of $G$ is contained in a unique
  maximal finitely generated Abelian subgroup, in particular, each
  Abelian subgroup of $G$ is finitely generated;
  \item \label{e:pr4} $G$ is finitely presented, and has only finitely many conjugacy classes
  of its maximal Abelian subgroups.
  \item \label{e:pr5} $G$ has solvable word problem, conjugacy problem and
  uniform membership problem.
  \item \label{e:pr7} $G$ has the Howson property. Namely, if $K_1$ and $K_2$
  are finitely generated subgroups of $G$, then the intersection $K_1\cap K_2$ is
  finitely generated. Moreover, for given finitely generated subgroups $K_1$ and $K_2$ of $G$,
  there is an algorithm to find the intersection $K_1\cap K_2$.
  \item \label{e:pr8} There is an algorithm to find the
  centralizer of a given element $g\in G$.
\end{enumerate}
\end{thm}
\begin{proof} Properties~\refe{pr0} and~\refe{pr1} follow immediately
from the definition of an $\mathcal{F}$-group. A proof of
property~\refe{pr3} can be found in~\cite{BMR};
property~\refe{pr2} is proven in~\cite{KM2}. Properties~\refe{pr2}
and ~\refe{pr4} are proved in~\cite{KM2}. Alternative proofs of
properties~\refe{pr3},~\refe{pr2} and~\refe{pr4} can be found
in~\cite{Sd}. Solvability of the word problem is shown
in~\cite{Mak84}, an algorithm to solve conjugacy problem can be
found in~\cite{MRS}. Recall that by a theorem proved by
Dahmani~\cite{Daco}, $\mathcal{F}$-groups are relatively
hyperbolic which allows one to use alternative algorithms to solve
word problem~\cite{Farb} and conjugacy problem~\cite{Bum}. Observe
that results proved in~\cite{Farb} imply finite presentability of
$\mathcal{F}$-groups, and a theorem proved in~\cite{Osin} implies
solvability of the conjugacy problem. Solvability of the uniform
membership problem and properties~\refe{pr7} and~\refe{pr8} are
proved in~\cite{KMRS}.
\end{proof}
The following Lemma~\ref{lem:elementary abelian splitting} asserts
that we can consider only those Abelian splittings of an
$\mathcal{F}$-group $G$ where each maximal Abelian non-cyclic
subgroup of $G$ is elliptic. We denote by $\mathcal{D}(G)$ the set
of all Abelian splittings of $G$ that have this latter property.
\begin{lem}\label{lem:elementary abelian splitting}
Let $G$ be an $\mathcal{F}$-group, let $M$ be a maximal Abelian
non-cyclic subgroup of $G$, and let $A$ be an Abelian subgroup of
$G$. If $G=G_1\ast_A G_2$, then $M$ can be conjugated into either
$G_1$ or $G_2$. If $G=G_1\ast_{A}$, and the intersection $M\cap
A^g$ is a proper subgroup of $M$ for some $g\in G$, then $M$ can
be conjugated so that $G=G_1\ast_{A}M$. If $G=G_1\ast_{A}$ and for
any $g\in G$, the intersection $M\cap A^g$ is either trivial or
coincides with $M$, then $M$ can be conjugated into $G_1$.
\end{lem}
\begin{proof} The first statement follows from the description of
commuting elements in a free product with amalgamation. Now, let
$G$ have the presentation as follows: $G=\langle G_v,t\mid
tat^{-1}=\omega(a)\forall a\in A\rangle$.

If $M\cap gAg^{-1}$ is not trivial, then by
Theorem~\ref{thm:properties}\refe{pr2}, $g^{-1}Mg$ is the maximal
Abelian subgroup containing $A$. Denote by $M_t$ the maximal
Abelian subgroup containing $tAt^{-1}$. Since the intersection
$g^{-1}Mg\cap t^{-1}M_tt=A$ is not trivial, by
Theorem~\ref{thm:properties}\refe{pr3}, we conclude that
$M_t=tg^{-1}Mgt^{-1}$. If $t\notin g^{-1}Mg$, then $A=g^{-1}Mg$,
so that $g^{-1}Mg$ is elliptic, as claimed. In this case,
$G=\langle G_v,t\mid tat^{-1}=\omega(a)\forall a\in M_1\rangle$
and $\omega(M_1)=M_2$ where both $M_1$ and $M_2$ are maximal
Abelian subgroups of $G_v$.

If $t\in g^{-1}Mg$, then $g^{-1}Mg\subseteq C_G(t)$, where
$C_G(t)$ is the centralizer of $t$ in $G$. According to the
presentation of $G$ as an HNN-extension, $C_G(t)=\langle
A,t\rangle\subseteq g^{-1}Mg$, hence $\langle
A,t\rangle=g^{-1}Mg$, in particular $A$ is a proper subgroup of
$g^{-1}Mg$ and $G=G_1\ast_A g^{-1}Mg$ (cf. also
\cite[Theorem~5]{GKM}).

If $M$ intersects no conjugate of $A$ and $M$ is hyperbolic when
acting on the Bass-Serre tree corresponding to the splitting of
$G$ as the HNN-extension, then $M$ inherits a non-trivial
splitting as a free product, a contradiction.
\end{proof}
\begin{defn}\label{defn:abelian cycle}
We say that an Abelian splitting
$\mathcal{S}=(\mathcal{G}(V,E);T,\theta)$ of a group $G$ is an
\emph{Abelian cycle of groups} if the following conditions hold:
\begin{enumerate}
  \item \label{e:ac3} $G$ can be obtained as a
  series of amalgamated products
  \[\tilde{G}=(((G_1\ast_{A_1}G_2)\ast_{A_2}G_3)\ast\dots)\ast_{A_{n-1}}G_n\]
  and an HNN-extension $G=\langle\tilde{G},t\mid
  A=t^{-1}\alpha(A_n)t\rangle$ with $A\subset G_1$ and $\alpha(A_n)\subset G_n$.
  In particular, the graph $(V,E)$ is a cycle.
  \item \label{e:ac4} The edge groups $A_1,\dots,A_n$ ($n\geq 1$) are all subgroups
  of a maximal Abelian subgroup $M\subset G$.
\end{enumerate}
We also call such a splitting $\mathcal{S}$ an $M$-\emph{cycle of
groups} to stress that all edge groups in $\Gamma$ are subgroups
of the group $M$.
 Thus, if $G$ is an $M$-cycle, then $G$ has the following
 presentation:
\begin{align*}
G=\langle G_1,\dots,G_n,t\mid \alpha(A_i)=\omega(A_i),
i=1,\dots,n-1, A=t^{-1}\alpha(A_n)t\rangle,
\end{align*}
where $\alpha(A_i)\subseteq G_i\cap M$ (for $i=1,\dots,n$),
$\omega(A_i)\subseteq G_{i+1}$ (for $i=1,\dots,n-1$) and $A\subset
G_1$.
\end{defn}
\begin{defn}\label{defn:star groups}
A graph of groups $\Upsilon(V,E)$ is a \emph{star of groups}, if
$(V,E)$ is a tree $T$ which has diameter $2$. If $\Upsilon$ is a
star of groups, then the fundamental group of $\Upsilon$ is as
follows:
\begin{align*}
\pi(\Upsilon)=\langle M,K_1,\dots,K_n\mid \alpha(A_i)=\omega(A_i),
i=1,\dots,n\rangle,
\end{align*}
meaning that $\alpha(A_i)\subseteq M$ and $\omega(A_i)\subseteq
K_i$. The vertex $v_0\in V$ with the stabilizer $M$ is called the
\emph{center} and vertices $v_i$ with stabilizers $K_i$ are called
\emph{leaves of} $\Upsilon(V,E)$. If $\Upsilon(V,E)$ is an Abelian
star of groups, then $M$ is a maximal Abelian subgroup of $G$.
\end{defn}
\begin{defn}\label{defn:constellation groups}
A graph of groups $\Psi(V,E\cup E_s)$ is a \emph{constellation of
groups}, if $\Psi(V,E\cup E_s)$ can be obtained by taking finitely
many amalgamated products of stars of groups over leaves and
HNN-extensions where both associated subgroups are stabilizers of
the centers of those stars. In other words, $\Psi(V,E\cup E_s)$
can be obtained by iterations of the following construction:
\[
\pi(\Psi)=\langle\pi(\Upsilon_1),\pi(\Upsilon_2),t\mid
K_i^{(1)}=K_j^{(2)}, M^{(1)}=tM^{(2)}t^{-1}\rangle,
\]
where $\pi(\Upsilon_l)=\langle
M^{(l)},K_1^{(l)},\dots,K_{n_l}^{(l)}\mid
\alpha(A_i^{(l)})=\omega(A_i^{(l)}), 1\leq i\leq n_l\rangle$ for
$l=1,2$ is a star of groups as in Definition~\ref{defn:star
groups}. We call an edge $e$ a \emph{silver edge} if $e$
corresponds to an HNN-extension where associated subgroups are
maximal Abelian. $E_s$ denotes the set of all silver edges in
$\Psi$.
\end{defn}
\begin{rem}\label{rem:constellation} In what follows, we focus on
Abelian stars of groups and constellations of groups, meaning that
edge groups are all Abelian.
\begin{enumerate}
 \item Since maximal Abelian subgroups of $G$ are malnormal by
 Theorem~\ref{thm:properties}\refe{pr3}, two Abelian stars of groups are
 never amalgamated over two different pairs of leaves, and the silver
 subgraph of $(V,E\cup E_s)$ is a tree.
 \item We \emph{do not} consider an amalgamated product of two stars
of groups with no HNN-extension a constellation of groups.
However, it is convenient to regard a star of groups as a
particular case of a (trivial) constellation of groups. We also
regard an edge of groups $M\ast_AG_v$ with $A\subseteq M$ and $M$
a maximal Abelian subgroup of $G$ as an Abelian star of groups.
\end{enumerate}
\end{rem}
\begin{lem}\label{lem:abelian cycle} If $G$ is an
$\mathcal{F}$-group and $\Delta(V,E)$ is a splitting of $G$ which
is an Abelian cycle of groups, then one can effectively modify
$\Delta$ so as to obtain a splitting $\Psi$ of $G$ which is an
Abelian constellation of groups.
\end{lem}
\begin{proof} Contract all edges of $\Delta$ but one to a point.
The new splitting of $G$ that we obtain is an HNN-extension
$G=G_v\ast_A$, hence $G$ has the presentation as follows:
$G=\langle G_v,t\mid tat^{-1}=\omega(a)\forall a\in A\rangle$. Let
$M$ be the maximal Abelian subgroup containing $A$.

First, assume that $A\varsubsetneqq M$. By
Lemma~\ref{lem:elementary abelian splitting}, $G=G_v\ast_AM$.
Furthermore, $G_v$ is an $\mathcal{F}$-group that splits into a
series of amalgamated products over Abelian subgroups. Observe
that all these Abelian subgroups and also $A$ are contained in a
maximal Abelian subgroup $M_v\subset G_v$.
Lemma~\ref{lem:elementary abelian splitting} implies that $M_v$
can be conjugated to a vertex group in the splitting of $G_v$, in
particular $A$ is elliptic in this splitting. Therefore, the
splitting of $G_v$ extends to a splitting of the whole group $G$
into a graph of groups that has a tree as the underlying graph,
with a vertex stabilized by $M$. Since all edge groups in the
graph are subgroups of $M$, by a sequence of slidings one obtains
a star of groups in the sense of Definition~\ref{defn:abelian
cycle}, as follows. If there is a vertex $v\in V$ such that
$M=G_v$, then define $u=v$, otherwise add a vertex $u$ with
$G_u=M$ and an edge $f$ with $G_f=M$ so that $i(f)=u$ and
$\tau(f)=v$ (this is a refinement of $\Delta$ at the vertex $v$).
Having introduced the vertex $u$ with the stabilizer $G_u=M$, we
make the following finite sequence of slidings in $\Delta$. Let
$v_i\in V$ be a vertex adjacent to $u$ (we set $v_1=v$), denote by
$f_i$ the edge connecting them (clearly, $f_1=f$), and assume that
$val(v_i)>1$ (for otherwise, we are done). Choose an edge $e\neq
f_i$ in $Star(v_i)$ and slide this edge to $u$. W.l.o.g., we can
assume that we had $i(e)=v_i$, so that having made the sliding we
have $i(e)=u$. If $G_{\tau(e)}\subseteq M$, then collapse $e$, so
that $u$ and $\tau(e)$ get identified. None of these operations
changes the fundamental group of $\Delta$. We end up with a star
of groups centered at $u$.

Now, let $A=M$. Since $M$ is malnormal in $G$, $M^t\neq M$ for
each $t\in G\setminus M$. Therefore, by the property~\refe{ac3} of
an Abelian cycle (see Definition~\ref{defn:abelian cycle} for the
notation), there is a unique edge $e\in E$ with $i(e)=v_n$ and
$\tau(e)=v_1$, so that the boundary monomorphisms are as follows:
$\alpha(G_e)=A_n=M$ and $\omega(G_e)=A=M^{t}$. To modify $\Delta$,
we add a vertex $u$ stabilized by $M$ and a vertex $u_t$
stabilized by $M^t$, join $u$ to $v_n$ by an edge $f_n$ with the
edge group $G_n=M$ and join $u_t$ to $v_1$ by an edge $f_t$ with
the edge group $G_t=M^t$. Next, we slide the edge $e$ along the
edges $f_n$ and $f_t$ so that $i(e)=u$ and $\tau(e)=u_t$; so $e$
becomes a \emph{silver} edge in the meaning of
Definition~\ref{defn:constellation groups}. Clearly, none of the
above operations changes the fundamental group of $\Delta$. The
graph spanned by the vertices $v_1,\dots,v_n,u$ is now a linear
tree (with no branch points) with all edge groups being subgroups
of $M$, hence one can transform this subgraph by a series of
slidings to an $M$-star of groups. Observe that $M\subset G_1$,
since $G_1$ contains $M^t$ and intersects with $M$ non-trivially.
Therefore, each edge group in this star of groups equals $M$. The
graph spanned by $v_1$ and $u_t$ is an edge of groups which is a
particular case of a star of groups with the center $u_t$ and a
unique leaf $v_1$. Thus, we have obtained a splitting $\Psi$ of
$G$ which is an Abelian constellation of groups.

It remains to notice that an Abelian $M$-cycle $\Delta$ can be
transformed to a constellation of groups (and not to a star of
groups) if and only if each edge group in $\Delta$ equals $M$.

To show that $\Psi$ can be found effectively, observe that we need
to use the following algorithms. First, for a given Abelian
subgroup $A$ of $G$ which is an edge group in a splitting of $G$,
one should find effectively the maximal Abelian subgroup $M$
containing $A$. Existence of this algorithm follows from
Theorem~\ref{thm:properties}~\refe{pr8}, as by
Theorem~\ref{thm:properties}~\refe{pr2}, $M$ is the centralizer of
any non-trivial element of $A$. The other problem which is to be
solved effectively is to find the intersection of two given
finitely generated subgroups of $G$. This algorithm is provided
according to Theorem~\ref{thm:properties}~\refe{pr7}.
\end{proof}
\begin{cor}\label{cor:two cycles} Let $G$ be an
$\mathcal{F}$-group, and let $M$ be a maximal Abelian subgroup of
$G$. If $G$ does not split as an HNN-extension where $M$ is one of
the two associated subgroups, then each splitting of $G$ contains
at most one Abelian $M$-cycle.
\end{cor}
\begin{proof} By the proof of Lemma~\ref{lem:abelian
cycle}, $G=G_v\ast_M$ if and only if $G$ has a splitting with an
Abelian $M$-cycle where $A_1=\dots=A_n=M$. Assume there are two
Abelian cycles in a splitting of $G$. One can find in each cycle
an edge (denoted by $e_1$ and $e_2$) that does not belong to the
other cycle, so that the edge group of both $e_1$ and $e_2$ are
proper subgroups of $M$. Choose a maximal tree $T$ in the
underlying graph so that $e_1,e_2\notin T$. Let $t_1$ and $t_2$ be
stable letters corresponding to $e_1$ and $e_2$. Since each edge
group in both cycles is a subgroup of $M$, according to the proof
of Lemma~\ref{lem:abelian cycle}, both $t_1$ and $t_2$ belong to
$M$, a contradiction.
\end{proof}
\subsection{Universal decomposition}
The following theorem~\ref{thm:effectiveJSJ} which is the main
result of~\cite{JSJ} is crucial for our proof. Before we state the
theorem, we need to introduce some more definitions.
\begin{defn}\label{defn:QHvertex}(QH-vertex)
Let $P$ be a planar subgroup of $G$ which admits one of the
following presentations:
\begin{enumerate}
  \item $\langle
  p_1,\dots,p_m,a_1,\dots,a_g,b_1,\dots,b_g\mid
  \prod_{k=1}^{m}p_k\prod_{j=1}^{g}[a_j,b_j]\rangle$;
  \item $\langle
  p_1,\dots,p_m,v_1,\dots,v_g\mid
  \prod_{k=1}^{m}p_k\prod_{j=1}^{g}v_j^2\rangle$.
\end{enumerate}
Let $\Gamma(V,E)$ be a graph of groups. Let $v\in V$ and let
$e_1,\dots,e_m$ be all edges with $i(e_i)=v$. We suppose that
$G_v=P$ and that $\alpha(e_i)=p_i$. Such a vertex $v$ is called a
\emph{QH-vertex}.
\end{defn}
\begin{defn}\label{defn:QHsbgp}(QH-subgroup) A subgroup $P$ of $G$ is a
\emph{QH-subgroup}, if there is a splitting $\mathcal{G}(V,E)$ of
$G$ and a QH-vertex $v\in\mathcal{G}$ (see
Definition~\ref{defn:QHvertex}) such that $P$ can be conjugated
into the stabilizer of $v$. A subgroup $P$ of $G$ is a
\emph{maximal QH-subgroup} (denoted by MQH-subgroup for short), if
for each elementary cyclic splitting $G=G_1\ast_C G_2$ either $P$
can be conjugated into $G_1$ or $G_2$, or $C$ can be conjugated
into $P$ in such a way that there is an elementary splitting of
$P$ over a cyclic subgroup $C_1$ so that this splitting extends to
an elementary splitting of the whole group $G$, and $C$ is
hyperbolic with respect to the splitting of $G$ over $C_1$.
\end{defn}
\begin{defn}\label{defn:almost reduced} We say that $\Delta$ is
\emph{almost reduced}, if the equality $G_e=G_v$
implies that $u$ is a QH-vertex (in particular, $G_e$ is cyclic),
$val(v)=2$ and for the other edge $f$ incident on $v$ we have that
$G_f\varsubsetneq G_v$ and the other endpoint of $f$ is a
QH-vertex as well.
\end{defn}
Recall that if $G$ is an $\mathcal{F}$-group, then
$\mathcal{D}(G)$ denotes the set of all Abelian splittings of $G$
where each maximal Abelian subgroup of $G$ is elliptic.
\begin{thm}\label{thm:effectiveJSJ}\cite[Theorem~0.1 and Proposition~2.15]{JSJ}.
Let $H$ be  a freely indecomposable $\mathcal{F}$-group. There
exists an almost reduced unfolded Abelian splitting $D\in{\mathcal
D}(H)$ of $H$ with the following properties:
\begin{enumerate}
  \item \label{e:effJSJ1} Every  MQH-subgroup of $H$ can be
  conjugated to a vertex group
in $D$; every QH-subgroup of $H$ can be conjugated into one of the
MQH-subgroups of $H$; non-MQH subgroups in $D$ are of two types:
maximal abelian and non-abelian, every non-MQH vertex group in $D$
is elliptic in every Abelian splitting in ${\mathcal D}(H)$.

   \item \label{e:effJSJ2} If an elementary cyclic splitting $H=A*_CB$ or
$H=A*_C$ is hyperbolic in another elementary cyclic splitting,
then $C$ can be conjugated into some MQH subgroup.

  \item \label{e:effJSJ3} Every elementary Abelian splitting $H=A*_CB$ or $H=A*_C$ from
${\mathcal D}(H)$ which is elliptic with respect to any other
elementary Abelian splitting from ${\mathcal D}(H)$ can be
obtained from $D$ by a sequence of collapses, foldings,
conjugations and modifying boundary monomorphisms by conjugation.

  \item \label{e:effJSJ4} If $D_1\in{\mathcal D}(H)$ is another splitting that has
properties~\refe{effJSJ1} and~\refe{effJSJ2}, then it can be
obtained from $D$ by slidings, conjugations, and modifying
boundary monomorphisms by conjugation.
\end{enumerate}
Moreover, given a presentation of $H$, there is an algorithm to
find the splitting $D$.
\end{thm}
In our proof, we use the slightly modified version of
Theorem~\ref{thm:effectiveJSJ}, stated in
Theorem~\ref{thm:universalJSJ} below. It follows
from~\cite[Theorem~6]{KM2} (cf. also~\cite[Theorem~4.1]{Sd}) that
an indecomposable $\mathcal{F}$-group $G$ is one of the following:
the fundamental group of a closed surface, a free Abelian or a
free non-Abelian group (cf. Example~\ref{ex:trivial decomp}).
\begin{thm}\label{thm:universalJSJ} Let $G$ be a one-ended
decomposable $\mathcal{F}$-group. $G$ has a semi-reduced Abelian
splitting $\Gamma=(\Gamma(V,E),T,\varphi)\in\mathcal{D}(G)$ called
a JSJ decomposition of $G$ that satisfies the following
properties:
\begin{enumerate}
  \item \label{e:i1} The decomposition $\Gamma$ is universal, in
  the meaning of Definition~\ref{defn:universal decomposition}.
  \item \label{e:i5} The Bass-Serre tree $\tilde{\Gamma}$
  corresponding to $\Gamma$ is $2$-acylindrical (see
  Definition~\ref{defn:Gtrees}).
  \item \label{e:i3} Each rigid vertex group in $\Gamma$ is of one of
the following two types: a maximal Abelian subgroup (we call such
a vertex elementary), or a non-Abelian subgroup.
  \item \label{e:i4} $(V,E)$ is a bipartite graph: two
  elementary vertices and two non-elementary vertices are never joined by an edge $e\in E$.
  \item \label{e:i2f} Each flexible vertex of $\Gamma$ is a
  maximal QH-vertex. Let $G=G_1\ast_C G_2$ or $G=G_1\ast_{C}$ be a
  cyclic splitting of $G$. $C$ can be conjugated into the stabilizer of a flexible
  vertex of $\Gamma$ if and only if the splitting in question is hyperbolic
  with respect to another splitting of $G$.
  \end{enumerate}
Moreover, there is an algorithm to obtain $\Gamma$.
\end{thm}
\begin{proof} Properties~\refe{i1} and~\refe{i2f} follow
immediately from Theorem~\ref{thm:effectiveJSJ} and the
definitions.

Let $\Delta\in\mathcal{D}(G)$ be an Abelian splitting which is the
output of the algorithm mentioned in
Theorem~\ref{thm:effectiveJSJ}. We modify the graph of groups
$\Delta$ so as to obtain a new splitting $\Gamma$ satisfying
properties~\refe{i3} and~\refe{i4}.

Let $M$ be a maximal Abelian subgroup of $G$ that contains either
$\alpha(G_e)$ or $\omega(G_e)$ for some $e\in E$. Consider the set
$E_M$ of all edges $e_i$ of $\Delta$ with $\alpha(A_i)\subseteq
M$. Since $M$ is elliptic in $\Delta$, the union $\Delta_M$ of all
edges $e\in E_M$ is a connected subgraph of $\Delta$. It is easy
to see that $\Delta_M$ can be found effectively. Indeed, it
follows from Theorem~\ref{thm:properties}~\refe{pr2} that an edge
with the stabilizer $A$ belongs to $\Delta_M$ if and only if a
non-trivial element of $A$ commutes with a non-trivial element of
$\alpha(A_1)$. By Theorem~\ref{thm:properties}~\refe{pr5}, the
word problem in $G$ is decidable so that this latter problem is
decidable also. If $\Delta_M$ is a tree, then by a series of
slidings it can be transformed to an $M$-star of groups (cf. the
proof of Lemma~\ref{lem:abelian cycle}). Otherwise, $\Delta_M$
contains Abelian cycles. It follows immediately from
Definition~\ref{defn:abelian cycle} that the union of all edges
$e_i$ of $\Delta_M$ with $\alpha(A_i)=\omega(A_i)$ form a maximal
tree of $\Delta_M$. Since $M^t\cap M^s=1$ for two different stable
letters $t\neq s$, the proof of Lemma~\ref{lem:abelian cycle}
shows that $\Delta_M$ can be transformed effectively into an
Abelian constellation of groups $\Psi_M$.

More generally, we have the following procedure. Since $M$ is
elliptic in $\Delta$, there exists a vertex $v\in V$ with
$M\subseteq G_v$. If $M\neq G_v$ for each $G_v$ that contains it,
then we add to $V$ an elementary vertex $z$ stabilized by $M$ and
connect $z$ by an edge $f$ with $G_f=M$ to $v$. When we have a
vertex for each maximal Abelian subgroup $M$, then we produce a
sequence of slidings as follows. If $\alpha(G_e)$ and
$\omega(G_e)$ are both subgroups of a maximal Abelian subgroup
$M$, then we slide $e$ so that $i(e)=z$ with $G_z=M$ and don't
change $\tau(e)$. If $\alpha(G_e)\subseteq M$ and
$\omega(G_e)\subseteq N$ for $N\neq M$, then we slide $e$ so that
$i(e)=z$ with $G_z=M$ and $i(e)=y$ with $G_y=N$ and declare $e$ a
\emph{silver edge}. The reason to introduce the more general
procedure is that in $\Gamma$ one can have cycles formed by an
$M$-tree and an $N$-tree. In this latter case we have silver edges
that do not belong to Abelian cycles in the sense of
Definition~\ref{defn:abelian cycle}. But the argument mentioned in
Remark~\ref{rem:constellation} remains valid in this case also,
and we conclude that the silver subgraph of the modified graph
$\Delta'$ is a forest. Therefore, we can collapse each silver
$M$-subtree to a point stabilized by $M$. Obviously, the
fundamental group of the new graph $\Gamma$ is isomorphic to $G$,
$\Gamma\in\mathcal{D}(G)$, and also properties~\refe{i3}
and~\refe{i4} hold. Furthermore, in $\Gamma$ each non-trivial
Abelian subgroup fixes a subgraph of diameter at most $2$ which,
together with the CSA property (see Theorem~\ref{thm:properties})
implies the assertion~\refe{i5}.
\end{proof}
\begin{cor}\label{cor:edge groups elliptic} Each edge group of
$\Gamma$ is elliptic in any splitting of $G$.
\end{cor}
\begin{proof} Let $\Lambda$ be a splitting of $G$, and let $G_e$ be an edge
stabilizer in $\Gamma$. We identify the edge $e$ with its lifting
to the Bass-Serre tree $\tilde{\Gamma}$. By
Theorem~\ref{thm:universalJSJ}~\refe{i1}, there is a
$G$-equivariant simplicial map $\kappa$ from $\tilde{\Gamma}$ onto
$\tilde{\Lambda}$. The image $\kappa(e)\in\Lambda$ of the edge
$e\in\tilde{\Gamma}$ is a path $\lambda$ in $\tilde{\Lambda}$; the
path $\lambda$ may be degenerate. As $\kappa$ is $G$-equivariant,
$G_e$ is a subgroup of the stabilizer $G_{\lambda}$ of $\lambda$,
in particular, $G_e$ fixes a point when acting on
$\tilde{\Lambda}$, hence is elliptic in $\Lambda$, as claimed.
\end{proof}
\begin{cor}\label{cor:basic properties JSJ} Let $T$ be
a simplicial $G$-tree so that $G$ acts on $T$ with Abelian edge
stabilizers.
\begin{enumerate}
 \item\label{e:cor1} Let $t\in T$ be an edge with the stabilizer $S_t$. If $S_t$
 is elliptic in any splitting of $G$, then $S_t$ can be conjugated into an
elementary vertex group of $\Gamma$.
 \item\label{e:cor2} If for each edge $t\in T$, the stabilizer $S_t$ is a subgroup of
$G$ which is elliptic in any splitting of $G$, then each flexible
vertex stabilizer of $\Gamma$ fixes a point in $T$.
\end{enumerate}
\end{cor}
\begin{proof} By Theorem~\ref{thm:universalJSJ}~\refe{i1}, there is a
$G$-equivariant simplicial map $\kappa$ from $\tilde{\Gamma}$ onto
$T$. If $e$ is an edge of $\tilde{\Gamma}$ such that $\kappa(e)$
contains $t$, then $G_e$ can be conjugated into $S_t$; in
particular, $G_e$ and a conjugate of $S_t$ belong to a maximal
Abelian subgroup of $G$. By
Theorem~\ref{thm:universalJSJ}~\refe{i4}, one of the two endpoints
of $e$ in $\Gamma$ is an elementary vertex with the stabilizer $M$
which is an Abelian subgroup of $G$. By Lemma~\ref{lem:abelian
cycle}, $M$ is a maximal Abelian subgroup of $G$, hence $S_t$ can
be conjugated into $M$, and the first assertion follows.

 To prove the second assertion, assume that a flexible vertex
stabilizer $G_u$ of $\Gamma$ does not fix a point in $T$. In this
case, $G_u$ inherits a non-trivial splitting $\Lambda$ from its
action on $T$. The edge groups in $\Lambda$ are subgroups of the
edge stabilizers of $T$. Collapse all the edges of $\Lambda$ but
one and denote by $\Lambda_1$ the obtained elementary splitting of
$G_u$. By Corollary~\ref{cor:edge groups elliptic}, the edge
groups of $G$ are elliptic when acting on $T$, so that $\Lambda_1$
extends to a splitting of $G$. Observe that the edge group of
$\Lambda_1$ is elliptic in any splitting of $G$, which contradicts
Theorem~\ref{thm:universalJSJ}\refe{i2f}.
\end{proof}
\subsection{Uniqueness of a universal decomposition}
\begin{lem} \label{lem:mapmu} Let $G$ and $H$ be two one-ended
$\mathcal{F}$-groups, and let
$\varphi\colon G\rightarrow H$ be an isomorphism. Let $\Gamma$ (or
$\Xi$) be an Abelian JSJ decomposition of $G$ (or $H$). Then there
exists a simplicial map $\mu\colon X\rightarrow Y$ between the
Bass-Serre trees $X=\tilde{\Gamma}$ and $Y=\tilde{\Xi}$ so that
the following diagram is commutative:
\[
\begin{CD}
{G\times X}        @> >>   {X}     \\
 @V{(\varphi,\mu)}VV                      @VV{\mu}V \\
{H\times Y}      @> >> {Y}
\end{CD}
\]
\end{lem}
\begin{proof} Observe that there are faithful actions $G\times Y\rightarrow Y$
defined by $\rho(g,y)=\varphi(g).y$ for all $g\in G$ and $y\in Y$,
and $H\times X\rightarrow X$ defined by
$\sigma(h,x)=\varphi^{-1}(h).x$ for all $h\in H$ and $x\in X$.
Furthermore, by Corollary~\ref{cor:edge groups elliptic}, each
edge group $H_e$ of $H$ fixes a point in $X$. Therefore, by
Corollary~\ref{cor:basic properties JSJ}\refe{cor1}, each flexible
vertex group of $X$ fixes a point (i.e., is elliptic) when acting
on $Y$. Observe that each rigid vertex group of $X$ is elliptic
also, by the definition. Each elementary vertex group $M$ of $X$
is a maximal Abelian subgroup of $G$, hence its image $\varphi(M)$
is a maximal Abelian subgroup of $H$. Since
$\Gamma\in\mathcal{D}(G)$ and $\Xi\in\mathcal{D}(H)$, $\varphi(M)$
fixes a vertex in $Y$. Moreover, since $G$ splits over a subgroup
$A\subseteq M$, we have that $H=\varphi(G)$ splits over
$\varphi(A)$, so that according to the proof of
Lemma~\ref{lem:abelian cycle} and Theorem~\ref{thm:universalJSJ},
$\varphi(M)$ fixes a unique elementary vertex in $Y$.

Our argument above allows one to define a simplicial map
$\mu\colon X\rightarrow Y$ as follows. If $v\in X$ is a vertex
with the stabilizer $G_v$, then $\mu(v)=y$ is the vertex with
$\varphi(G_v)\subseteq H_y$. If $e\in X$ is an edge with the
endpoints $v$ and $u$, then $\mu(e)=f$ is the path joining
$\mu(v)$ and $\mu(u)$. Furthermore, we claim that the diagram in
the assertion of the theorem is commutative. Let $g\in G$ be a
non-trivial element, and let $v\in X$ be a vertex with the
stabilizer $G_v$. The image $u=g.v\in X$ is the vertex with the
stabilizer $G_u=g^{-1}G_vg$, hence $\mu(g.v)=y_u\in Y$ so that
$\varphi(g^{-1}G_vg)\subseteq H_u$, where $H_u$ denotes the
stabilizer of $y_u$. On the other hand, $\mu(v)=y_v$ with
$\varphi(G_v)\subseteq H_v$, and $g$ maps $y_v\in Y$ to the vertex
$\bar{y}_v=\varphi(g).y_v$ with the stabilizer
$H_{\bar{v}}=\varphi(g)^{-1}H_v\varphi(g)$. Observe that both
$H_u$ and $H_{\bar{v}}$ contain $\varphi(g^{-1}G_vg)$ as a
subgroup. If $G_v$ (hence, $\varphi(G_v)$) is non-elementary, then
it cannot fix an edge in either $X$ or $Y$. If $G_v$ is
elementary, then $g.v$, $y_u=\mu(g.v)$, $y_v=\mu(v)$ and
$\varphi(g).\mu(v)$ are elementary vertices. In either case, we
conclude that $H_u=H_{\bar{v}}$, and since the vertex of $Y$
stabilized by $H_u$ is unique, we have that
$\mu(g.v)=\varphi(g).\mu(v)$, as claimed.
\end{proof}
\begin{thm} \label{thm:down_to_vertices} Let $\varphi\colon G\rightarrow
H$ be an isomorphism of two one-ended $\mathcal{F}$-groups, and
let $\Gamma=\Gamma(V,E)$ and $\Xi=\Xi(U,B)$ be Abelian JSJ
decompositions of $G$ and $H$, respectively. Then the equivariant
map $\mu\colon\tilde{\Gamma}\rightarrow\tilde{\Xi}$ between the
Bass-Serre trees, defined in Lemma~\ref{lem:mapmu}, is a
one-to-one isometry.
\end{thm}
\begin{proof} Denote $X=\tilde{\Gamma}$ and $Y=\tilde{\Xi}$.
First, observe that the length of the image $\mu(e)\in Y$ of an
edge $e\in X$ does not exceed $2$ since $Y$ is $2$-acylindrical.
Moreover, according to Theorem~\ref{thm:universalJSJ}~\refe{i4},
we can assume that one of the endpoints $u$ and $v$ of $e$ is an
elementary vertex, so that the image of this endpoint in $Y$ is an
elementary vertex as well. As $\Xi$ is a bipartite graph (hence,
$Y$ is a bipartite tree) and different elementary vertex
stabilizers have only trivial intersections, it follows that
$\mu(e)$ has length $1$ or $0$.

Now, we claim that the non-degenerate images of two edges of $X$
cannot get folded in $Y$. More precisely, let $e$ and $f$ be two
edges of $X$, both incident on a vertex $v$ so that $i(e)=i(f)=v$,
hence $\alpha(G_e),\alpha(G_f)\subseteq G_v$, with different
terminal points: $\tau(e)=u$ and $\tau(f)=w$. Assume that the
images of $e$ and $f$ under $\mu$ get folded, so that
$\mu(e)=\mu(f)=c$ and $\mu(u)=\mu(w)=y$. Let $\mu(v)=y_v$. Since
$\varphi(g).\mu(x)=\mu(g.x)$ for all $x\in X$ and $g\in G$, both
$\varphi(G_e)$ and $\varphi(G_f)$ are subgroups of $H_c$. Since
the edge stabilizers in $Y$ are Abelian, $H_c$ is an Abelian
subgroup of $H$ and therefore, is a subgroup of a unique maximal
Abelian subgroup of $H$ which we denote by $M_H$. It follows that
both $\varphi(G_e)$ and $\varphi(G_f)$ are subgroups of $M_H$, so
that both $G_e$ and $G_f$ are subgroups of a maximal Abelian
subgroup $M=\varphi^{-1}(M_H)$. Hence, by our construction, $v$ is
an elementary vertex of $X$ (and $M=G_v$). Therefore, $G_u$ and
$G_w$ are non-elementary, and neither is $H_y$ as both
$\varphi(G_u)$ and $\varphi(G_w)$ are subgroups of $H_y$. On the
other hand, $H_y$ inherits a non-trivial elementary splitting from
its action on $X$, a contradiction.

Next, we show that the image of an edge $e\in X$ cannot have
length $0$ in $Y$. Assume that $\mu(v)=\mu(u)=y$, where $u$ and
$v$ are the endpoints of $e$, and $i(e)=v$ is an elementary
vertex. If $\alpha(G_e)\varsubsetneqq G_v$, then we get a
contradiction, because $H_y$ acts non-trivially on $X$, hence
splits over an Abelian subgroup. Let $G_e=G_v$. In this case the
valence of $v$ is at least $2$ since $X$ is semi-reduced; let
$f\neq e$ be another edge incident on $v$. As we have just shown,
the images of edges incident on an elementary vertex in $X$ cannot
get folded in $Y$. If the image of $f$ under $\mu$ collapses also,
then we have three vertices of $X$ mapped to a vertex $y\in Y$, so
that $H_y$ acts non-trivially on $X$, a contradiction. Thus,
$\mu(f)$ is not degenerate, so that in $Y$ there is an edge
stabilized by $\varphi(G_f)\subset \varphi(G_v)$. By our
construction of the graph $\Xi$ in Theorem~\ref{thm:universalJSJ},
there is an elementary vertex $z$ in $Y$ with the stabilizer
$\varphi(G_v)$. By the definition of $\mu$, $z=\mu(v)\neq \mu(u)$,
a contradiction.

So far, we have shown that $\mu$ is a local immersion. Finally,
assume that there are two edges (or vertices) of $X$ which are
mapped to the same edge (or vertex) in $Y$. Consider the path $p$
connecting them in $X$ and its image $\mu(p)$ in $Y$. Since
$\mu(p)$ is a closed path in $Y$, $p$ has either an edge $e$
incident on a vertex $v$ so that $\mu(e)=\mu(v)$, or two edges $e$
and $f$ incident on $v$ so that $\mu(e)=\mu(f)$. In either case,
$\mu$ restricted to $Star(v)$ is not a local immersion, a
contradiction.
\end{proof}
\subsection{Isomorphism of groups and splittings of the groups}\label{s:iso split}
\begin{thm}\label{thm:isomorphic image of JSJ} Let $G$ and $H$ be two
$\mathcal{F}$-groups, and let $\varphi\colon G\rightarrow H$ be an
isomorphism. Let $\Gamma=\Gamma(V,E)$ (or $\Xi=\Xi(U,B)$) be the
Abelian JSJ decomposition of $G$ (or $H$). The image
$\varphi_*(\Gamma)$ of $\Gamma$ under $\varphi$ can be obtained
from $\Xi$ by conjugation and modifying boundary monomorphisms.
\end{thm}
\begin{proof} Fix the natural lift $D$ of $\Gamma$ into
$X=\tilde{\Gamma}$ (see Definition~\ref{defn:natural lift}), and
let $\mu(D)$ be the image of $D$ in $Y=\tilde{\Xi}$, where $\mu$
is the $G$-equivariant isometry defined in Lemma~\ref{lem:mapmu}
(see also Theorem~\ref{thm:down_to_vertices}); recall that $G$
acts on $X$ by left multiplications and on $Y$ via the isomorphism
$\varphi$ and left multiplications. Observe that $\mu(D)$ is a
fundamental domain of $Y$. Indeed, since $X=G.D$, and the map
$\mu$ is $G$-equivariant and onto, we conclude that $Y=G.\mu(D)$.
Moreover, as $\mu$ is $G$-equivariant, $x_1=g.x_2$ iff
$\mu(x_1)=\varphi(g).\mu(x_2)$ for all $x_1,x_2\in X$ and $g\in
G$, so that two vertices (or two edges) of $X$ are $G$-equivalent
if and only if their images in $Y$ are $\varphi(G)=H$-equivalent.
Therefore, two different edges of $\mu(D)$ are never
$H$-equivalent, and two vertices $\mu(v)$ and $\mu(u)$ of $\mu(D)$
are $H$-equivalent if and only if $v$ and $u$ are $G$-equivalent.
This latter argument shows that the underlying graph of $\Gamma$
is the underlying graph of $\Xi$. Therefore, we can assume that
the maximal trees of $\Gamma$ and of $\Xi$ coincide and the
orientation of edges is the same. It can be readily seen that
$\varphi_*(\Gamma)$ can be obtained from $\mu(D)$ by identifying
the $H$-equivalent vertices.

Now, let $K$ be the natural lift of $\Xi$ into $Y$. Fix a vertex
$d\in D$. There is a vertex $k\in K$ so that $d$ and $k$ are
$G$-equivalent. Observe (cf. Lemma~\ref{lem:ext domain}) that
there may be more than one vertex $G$-equivalent to $d$. To
specify our choice, we also require that there is an isomorphism
of graphs $\lambda$ with $\lambda(k)=d$ that maps each vertex (or
edge) of $K$ to a $G$-equivalent vertex (or edge) in $D$.

The stabilizers $H_d$ and $H_k$ are conjugate in $H$, let
$H_d=H_k^h$ for some $h\in H$. Let $e_k$ be an edge in $K$
incident on $k$, and let $k_1$ be the other endpoint of $e_k$.
Recall that by our construction, precisely one of the vertices $k$
and $k_1$ is elementary, so that $k$ and $k_1$ are never
$G$-equivalent. Denote $e_d=\lambda(e_k)$ and $d_1=\lambda(k_1)$.
We have that $H_{d_1}=H_{k_1}^{h_1}$ for some $h_1\in H$, so that
$H_{e_d}=H_d\cap H_{d_1}=H_k^h\cap H_{k_1}^{h_1}=(H_k\cap
H_{k_1}^{h_1h^{-1}})^h=(H_k^{hh_1^{-1}}\cap H_{k_1})^{h_1}$. Since
the tree $Y$ is $2$-acylindrical and $\lambda$ is an isometry,
this latter intersection is non-empty if and only if either
$h_1h^{-1}\in H_k$ so that
$H_{e_d}=(H_{e_k}^{h_1h^{-1}})^h=H_{e_k}^{h_1}$, or $h_1h^{-1}\in
H_{k_1}$ so that $H_{e_d}=H_{e_k}^{h}$. In either case, we need to
modify a boundary monomorphism.

Consider a particular case when the natural projections of $d$ and
$d_1$ into $\Gamma$ are joined by two edges. We use  the above
notation. Let $f\neq e_d$ be the other edge joining $d$ and $d_1$
in $(V,E)$, and let $t$ be the stable letter that corresponds to
$f$ in $\Gamma$. W.l.o.g., we can assume that $i(f)=d$. Since the
graphs $(V,E)$ and $(U,P)$ are isomorphic, there is a unique edge
$p\in P$ so that $\lambda(f)=p$: this is the edge joining (the
natural projections of) $k$ and $k_1$ in $(U,P)$. We denote by $s$
the stable letter that corresponds to $p$ in $\Xi$. Let
$A=\alpha(G_f)\subseteq G_d$ and $B=\omega(G_f)\subseteq G_{d_1}$,
so that $A^t=B$. As we have just shown, $l=h_1h^{-1}$ is either in
$H_k$ or in $H_{k_1}$. If $l\in H_k$, then $s=\varphi(t)lh_f$ with
$h_f\in H_k$ being non-trivial if we need to modify the boundary
monomorphism as follows:
$\alpha(H_{\mu(f)})=h_f\alpha(H_p)h_f^{-1}$. If $l\in H_{k_1}$,
then $s=l^{-1}\varphi(t)h_f$ with $h_f\in H_k^l$ and
$\alpha(H_{\mu(f)})=h_f\alpha(H_p)^lh_f^{-1}$.

We proceed with the other edges incident on $k$ and check that the
assertion holds for $Star(k)$. The assertion follows by induction
on the number of vertices.
\end{proof}

\section{Algorithm to solve the isomorphism problem}\label{s:algorithm}
Our algorithm is based on the following result.
\begin{thm}\cite[Theorem~0.1 and Theorem~13.1]{JSJ}\label{thm:MRmachine} Let $\langle
\mathcal{S}\mid\mathcal{R}\rangle$ be a finite presentation of an
$\mathcal{F}$-group $G$; we regard this presentation as the input
of Elimination process. The Elimination process determines whether
or not $G$ is freely indecomposable, and the output of the process
is a finite presentation $\langle S\mid R\rangle$ of $G$ that can
be described as follows:
\begin{enumerate}
 \item\label{e:MRfree} If $G$ is a free non-Abelian group, then
 $R=\emptyset$.
 \item\label{e:MRfreeprod} If $G$ is freely decomposable but not free,
 then there are partitions $S=S_1\sqcup...\sqcup S_k\sqcup S_{k+1}$ and $R=R_1\sqcup...\sqcup
 R_k$, so that $\langle S\mid R\rangle=\langle S_1\mid R_1\rangle\ast\dots\ast\langle S_k\mid
 R_k\rangle\ast\langle S_{k+1}\mid -\rangle$, where $\langle S_i\mid R_i\rangle$ is a presentation
 of a freely indecomposable non-cyclic group for $1\leq i\leq k$, and
 $\# S_{k+1}\geq 0$. In other words, the presentation $\langle S\mid
 R\rangle$ corresponds to the Grushko decomposition of $G$.
 \item\label{e:MRJSJ} If $G$ is freely indecomposable, then the
 output of the Elimination process is a presentation of $G$ as a JSJ-graph
 of groups. If $G$ is also indecomposable in the meaning of
 Definition~\ref{defn:indecomposable group}, then the presentation
 $\langle S\mid R\rangle$ of $G$ has the following properties.
  \begin{enumerate}
  \item\label{e:MRcs} If $G$ is the fundamental group of a closed
 surface, then $R$ is a set of quadratic words, in the standard form.
  \item\label{e:MRab} If $G$ is a free Abelian group, then the
 cardinality of $S$ is minimum possible; in other words, $\#
 S=rank(G)$.
  \end{enumerate}
\end{enumerate}
\end{thm}
In what follows, we assume that we are given a presentation of
$G=\langle \mathcal{S}_G\mid\mathcal{R}_G\rangle$ and a
presentation of $H=\langle \mathcal{S}_H\mid\mathcal{R}_H\rangle$,
both presentations are output of the Elimination process.

\begin{lem} \label{lem:trivialJSJ} Let $G$ and $H$ be
indecomposable $\mathcal{F}$-groups. There exists an effective
procedure to decide whether or not $G$ and $H$ are isomorphic.
\end{lem}
\begin{proof} We apply the Elimination process to both presentations of $G$
and of $H$ to determine whether or not the corresponding group is
a free group. If both $G$ and $H$ are free, then they are
isomorphic if and only if the cardinalities of their generating
sets coincide. Now, assume that neither of $G$ and $H$ is a free
group. Since the equalities $[g_i,g_j]=1$ for all pairs of
generators of $G$ hold in $G$ if and only if $G$ is a free Abelian
group, and the word problem for $\mathcal{F}$-groups is solvable
by Theorem~\ref{thm:properties}\refe{pr5}, one can effectively
decide whether or not $G$ and $H$ are free Abelian groups.
Moreover, if $G$ is a free Abelian group, then by
Theorem~\ref{thm:MRmachine}\refe{MRab}, one can effectively
determine the rank of $G$. If both groups $G$ and $H$ are free
Abelian, then they are isomorphic if and only if their ranks are
equal. If neither of $G$ and $H$ is free Abelian, then both $G$
and $H$ are fundamental groups of closed surfaces. By
Theorem~\ref{thm:MRmachine}\refe{MRcs}, one can effectively find
standard quadratic presentations for both $G$ and $H$. The groups
are isomorphic if and only if their standard presentations
coincide, up to permutation of generators.
\end{proof}
In what follows, we assume that both $G$ and $H$ are decomposable
groups. Lemma~\ref{lem:free decomposition} below allows us to
reduce the problem to the case when both $G$ and $H$ are freely
indecomposable groups.
\begin{lem}\label{lem:free decomposition}\cite{Ku}
Let $G=G_1\ast G_2\ast...\ast G_k\ast F_r$ and $H=H_1\ast
H_2\ast...\ast H_l\ast F_s$ be the Grushko decompositions. The
groups $G$ and $H$ are isomorphic if and only if $k=l$, $r=s$ and
there exists a permutation $\sigma$ of the set $\{1,\dots,k\}$ so
that $G_i$ is isomorphic to $H_{\sigma(i)}$ for each
$i=1,\dots,k$.
\end{lem}

\subsection{Freely indecomposable groups}  Our solution to the
isomorphism problem relies upon Theorem~\ref{thm:isomorphic image
of JSJ}.  According to Theorem~\ref{thm:MRmachine}\refe{MRJSJ},
the above presentations define $G$ and $H$ as fundamental groups
of graphs of groups: $G\backsimeq\pi_1(\Gamma)$ and
$H\backsimeq\pi_1(\Xi)$, which are Abelian JSJ decompositions of
$G$ and $H$, respectively. Our algorithm is built so as to compare
the two graphs of groups and conclude whether or not their
fundamental groups are isomorphic; the algorithm is described in
Theorem~\ref{thm:algorithm} below. It consists of a sequence of
smaller procedures, some of these we describe now. First, we
classify and compare the vertex groups.
\begin{lem}\label{lem:classify} There is an algorithm to determine
the type of a given vertex $G_v$ in an Abelian JSJ decomposition
$\Gamma$ of an $\mathcal{F}$-group.
\end{lem}
\begin{proof} If each pair of generators commute, then $G_v$ is
free Abelian. If $G_v$ is flexible, then the given presentation of
$G_v$ is a presentation of a QH-subgroup of one of the two
possible kinds~\ref{defn:QHsbgp}, up to permutation of the
generators. If $G_v$ is neither Abelian nor flexible, then
according to Theorem~\ref{thm:universalJSJ}, $G_v$ is rigid
non-elementary.
\end{proof}
\begin{defn}\label{defn:extendable iso} Let $G$ and $H$ be two
isomorphic groups, let $A_1,\dots,A_n$ be subgroups of $G$, and
let $B_1,\dots,B_n$ be subgroups of $H$. An isomorphism
$\phi\colon G\rightarrow H$ is an \emph{extendable isomorphism}
(or \emph{e-isomorphism }for short), if there is one-to-one
correspondence $A_i\rightarrow B_{j_i}$ between the sets of the
subgroups so that $\phi$ maps $A_i$ onto a conjugate of $B_{j_i}$.
Pairs $(G,\{A_1,\dots,A_n\})$ and $(H,\{B_1,\dots,B_n\})$ are
called \emph{e-isomorphic}, if there is an e-isomorphism
$\phi\colon G\rightarrow H$.
\end{defn}

%
To find e-isomorphisms of QH-subgroups, we use the Elimination
process that gives their standard presentations, and the following
classical result.
\begin{lem}\label{lem:flexible gps} Let $G_v\subset G$ and
$H_u\subset H$ be two QH-subgroups in the Abelian JSJ
decompositions of one-ended $\mathcal{F}$-groups $G$ and $H$, and
let $A_1,\dots,A_n\subset G_v$ and $B_1,\dots,B_n\subset G_u$ be
their sets of peripheral subgroups. Then $G_v$ and $H_u$ are
e-isomorphic if and only if their standard presentations (see
Definition~\ref{defn:QHsbgp}) are the same, up to permutation of
generators. In particular, if $\varphi_v$ is an e-isomorphism,
then $\varphi_v(A_i)=B_i$ for all $i=1,2,\dots,n$.
\end{lem}
%
\subsection{Rigid vertices}
 To find out whether or not two rigid
vertex groups are e-isomorphic, we use Theorem~\ref{thm:monom}
below. To state the theorem, we need some more definitions.
\begin{defn}\label{defn:equivalent mono} Two monomorphisms $\psi\colon
G\rightarrow H$ and $\phi\colon G\rightarrow H$ are
\emph{equivalent} if $\psi$ is a composition of $\phi$ and
conjugation by an element from $H$.
\end{defn}
\begin{defn}\label{defn:splitting modulo} Let $G$ be a group
and ${\mathcal K} = \{K_1,\ldots,K_n\}$ be a set of subgroups of
$G$. An Abelian splitting $\Delta$ of $G$ is called a {\em
splitting modulo} $\mathcal K$ if all subgroups from ${\mathcal
K}$ are conjugated into vertex groups in $\Delta$.
\end{defn}
Observe that a rigid vertex group in an Abelian JSJ  decomposition
of a group has no non-degenerate Abelian splittings modulo its
peripheral subgroups.
\begin{thm} \cite[Theorem~15.1]{JSJ}\label{thm:monom} Let $G$ (or $H$) be an $\mathcal{F}$-group,
and let $S_A=\{A_1,\ldots,A_n\}$ (respectively,
$S_B=\{B_1,\ldots,B_n\}$) be a finite set of non-conjugated
maximal Abelian subgroups of $G$ (respectively, $H$) such that the
Abelian decomposition of $G$ modulo $S_A$ is trivial. The number
of equivalence classes of monomorphisms from $G$ to $H$ that map
subgroups from $S_A$ onto  conjugates of the corresponding
subgroups from $S_B$ is finite. A set of representatives of the
equivalence classes can be effectively found.\end{thm}
\begin{cor}\label{cor:BPproperty} Let $G$ be an
$\mathcal{F}$-group, and let $S=\{A_1,\ldots ,A_n\}$ be a finite
set of maximal Abelian subgroups of $G$. Denote by $Out(G;S)$ the
set of those outer automorphisms of $G$ which map each $A_i\in S$
onto a conjugate of it. If $Out(G;S)$ is infinite, then $G$ has a
non-trivial Abelian splitting modulo $S$. There is an algorithm to
decide if $Out(G;S)$ is infinite and if it is, to find the
splitting.
\end{cor}
\begin{lem}\label{lem:rigid gps} Let $G$ (or $H$) be an $\mathcal{F}$-group,
and let $S_A=\{A_1,\ldots,A_n\}$ (respectively,
$S_B=\{B_1,\ldots,B_n\}$) be a finite set of non-conjugated
maximal Abelian subgroups of $G$ (respectively, $H$) such that the
Abelian decomposition of $G$ modulo $S_A$ is trivial. Then there
is an algorithm to decide whether or not $G$ and $H$ are
e-isomorphic, and if they are, then the algorithm finds all the
equivalence classes of extendable isomorphisms from $G$ to $H$.
\end{lem}
\begin{proof} We apply Theorem \ref{thm:monom} and find all the representatives
$\phi _1,\ldots ,\phi _k$ (if exist) of the equivalence classes of
monomorphisms from $G_v$ to $H_u$ that map subgroups from $S_A$
onto the subgroups from $S_B$.

If a monomorphism $\phi\colon G_v\rightarrow H_u$ that maps the
edge groups of $G_v$ onto the conjugates of the corresponding edge
groups of $H_u$ exists, one can effectively check whether or not
it is onto. First, we apply~\cite[Theorem 3.21]{JSJ} to obtain a
presentation for the image $\phi(G)\subseteq H$ which is the
subgroup of $H$ generated by $x_1,\dots,x_k$. Now, we
apply~\cite{MRS}  to see whether or not $h_j\in\phi(G)$ for each
$j$. The monomorphism $\phi$ is onto if and only if  $\phi$ is an
isomorphism.
\end{proof}
%
\subsection{Algorithm}\label{s:alg}
Let $\hat{\Gamma}(V,E)$ and $\hat{\Xi}(U,P)$ be Abelian JSJ
decompositions of two one-ended $\mathcal{F}$-groups $G$ and $H$,
respectively (see Theorem~\ref{thm:universalJSJ}). Assume that
there is an isomorphism of graphs $\lambda\colon (V,E)\rightarrow
(U,P)$. We denote the image $\lambda(e)$ of an edge $e$ by the
same letter $e$. For each vertex $v\in V$, we order all the edges
incident on $v$ end fix the same order for the edge subgroups of
$G_v$, so that $A_i=\alpha(G_{e_i})$. (Since our ordering is local
and the graph is bipartite, we can always assume that $i(e_i)=v$.)
Similarly, we order all the edge subgroups of $H_u$ where
$u=\lambda(v)$ when we assume that $\lambda$ respects the ordering
of edges incident on $v$ and on $u$. Further, we assume that for
each $v\in V$ and $u=\lambda(v)$, there is an e-isomorphism
$\varphi_v\colon G_v\rightarrow H_{u}$ that preserves ordering of
the edge subgroups of $G_v$ and $G_u$, so that $\varphi_v(A_i)$ is
conjugate to $B_i$ in $H_u$.

We fix a maximal tree $T$ in $(V,E)$ (hence, in $(U,P)$) and
introduce comparative labelling of edges $L^{(\varphi)}_u\colon
P\cap T\rightarrow H_u$ defined as follows. Let $v\in V$ be a
\emph{rigid non-elementary} vertex, and let $A_1,\dots, A_n\subset
G_v$ be the edge subgroups. For $u=\lambda(v)$, let $B_1,\dots,
B_n\subset H_u$ be the edge subgroups. Fix an e-isomorphism
$\varphi\colon G_v\rightarrow H_u$ and set
$L^{(\varphi)}_u(p_i)=h_i\in H$ if $p_i\in P$ is an edge incident
on $u$ with the edge group $B_i$ and
$\varphi(A_i)=h_iB_ih_i^{-1}$. Notice that labelling depends on
the e-isomorphism $\varphi$. We assign the trivial label $1\in H$
to each edge $e\in T$ incident on a flexible vertex.
By a \emph{star of a vertex} $v$ in the tree $T$ we mean the
subgraph $Star(v)$ of $T$ where the set of edges consists of the
edges of $T$ incident on $v$ and the set of vertices consists of
the endpoints of those edges.
\begin{lem}\label{lem:criterion} With the above notation and assumptions,
e-isomorphisms between vertices of $\Gamma(V,E)$ and $\Xi(U,B)$
can be extended to an isomorphism between the fundamental groups
$\pi_1(\Gamma)$ and $\pi_1(\Xi)$ if and only if there are
e-isomorphisms of vertices so that in the star of each elementary
vertex, at most one label $h_i$ is not trivial.
\end{lem}
\begin{proof} To show that the condition is necessary, suppose there is
an elementary vertex $u$ with two different edges $e_1,e_2\in
Star(u)$ stabilized by $B_1,B_2$, so that their labels $h_1$ and
$h_2$ are not trivial. Observe that $B_1,B_2\subset H_u$, so that
$B_i^{h_i}\subset H_u^{h_i}$ for $i=1,2$. Therefore,
$H_u^{h_1}=H_u^{h_2}$, hence $h_1h_2^{-1}\in H_u$, a
contradiction.

To show that the condition is also sufficient, we extend
e-isomorphisms $\varphi_v\colon G_v\rightarrow H_{\lambda(v)}$
between vertices of the graphs of groups $\Gamma(V,E)$ and
$\Xi(U,B)$ to an isomorphism between the fundamental groups of the
trees of groups
$\varphi_T\colon\hat{\Gamma}(T)\rightarrow\hat{\Xi}(T)$. These
trees of groups are obtained from the graphs of groups $\Gamma$
and $\Xi$ by removing the edges that do not belong to $T$. The map
$\varphi_T$ defines the images of the vertex groups $G_v$ of
$\Gamma$ under an isomorphism $\varphi\colon G\rightarrow H$ that
we are constructing. Having defined images of $G_v$ in $H$, we
assign images to the stable letters in the presentation of $G$ as
the fundamental group of $\Gamma(V,E)$, and get the isomorphism
$\varphi\colon G\rightarrow H$.

Fix elementary vertices $u\in U$ and $v\in V$ so that
$u=\lambda(v)$.
%
First, we extend e-isomorphisms between vertices of $Star(v)$ and
$Star(u)$ to an e-isomorphism between the fundamental groups
$\pi_1(Star(v))$ and $\pi_1(Star(u))$. Assume that in $Star(u)$,
precisely one label $h$ is not trivial. Let $u_o=\tau(e_0)$ where
$e_0$ is the labelled edge, and let $T_u$ denote the connected
component of $T\setminus\{e_0\}$ that contains $u$. We replace the
e-isomorphism $\varphi_x\colon G_x\rightarrow H_{\lambda(x)}$ by
$\hat{h}\circ\varphi_x$ where $\hat{h}$ is conjugation by $h$, for
each $x$ with $\lambda(x)\in T_u$. Let $v_0\in V$ be so that
$u_0=\lambda(v_0)$ and $\varphi_v^{(0)}\colon G_{v_0}\rightarrow
H_{u_0}$ be the e-isomorphism that corresponds to the labelling in
question. Observe that all vertices of $Star(u)$ but $u_0$ are in
$T_u$, and e-isomorphisms $\hat{h}\circ\varphi_v$ and
$\phi\in\{\varphi_v^{(0)},\hat{h}\circ\varphi_x\mid x\in Star(v),
x\neq v_0\}$ agree on edge subgroups.
Therefore, the e-isomorphisms $\hat{h}\circ\varphi_v$ and
$\varphi_v^{(0)},\hat{h}\circ\varphi_x\ (x\in Star(v), x\neq v_0)$
define an e-isomorphism $\psi_v$ between the fundamental groups
$\pi_1(Star(v))$ and $\pi_1(Star(u))$, since replacing
e-isomorphisms at the vertices $x\in V$ with
$\lambda(x)\in\Delta_u$, does not affect the labelling of $P$. If
there is no non-trivial label in $Star(u)$, then the
e-isomorphisms $\varphi_x$ where $x\in Star(v)$ agree on edge
subgroups, hence extend to an e-isomorphism
$\psi_v\colon\pi_1(Star(v))\rightarrow\pi_1(Star(u))$.

We proceed to other elementary vertices by induction on the
distance from $v$ in $T$ and end up with the isomorphism
$\varphi_T$. Now, let $e\in E$ do not belong to $T$, and let $t$
(or $s$) be the stable letter that corresponds to $e$ in $G$ (or
$H$). Let $v=i(e)$ and $x=\tau(e)$ be the endpoints of $e$, and
$A=\alpha(G_e)$ and $C=\omega(G_e)$. Recall that
$\alpha(G_e)=\omega(G_e)^t$ in $G$ and $\alpha(H_e)=\omega(H_e)^s$
in $H$. Our assumptions and the above procedure imply that
$\varphi_T(\alpha(G_e))=\alpha(H_e)^h$ and
$\varphi_T(\omega(G_e))=\omega(H_e)^b$ for some $h$ and $b$ in
$H$. Hence, we can set $\varphi_T(t)=b^{-1}sh$ to preserve the
relations. Obviously, the map $\varphi\colon G\rightarrow H$ that
we obtain is an isomorphism.
\end{proof}
\begin{thm}\label{thm:algorithm} Let $G\cong\langle
\mathcal{S}_G\mid\mathcal{R}_G\rangle$ and $H\cong\langle
\mathcal{S}_H\mid\mathcal{R}_H\rangle$ be finite presentations of
fully residually free groups. There exists an algorithm that
determines whether or not $G$ and $H$ are isomorphic. If the
groups are isomorphic, then the algorithm finds an isomorphism
$G\rightarrow H$.
\end{thm}
\begin{proof} We apply the Elimination process to the given presentations.
The output of the Elimination process are presentations
$G\cong\langle S_G\mid R_G\rangle$ and $H\cong\langle S_H\mid
R_H\rangle$ described in Theorem~\ref{thm:MRmachine}. If both $G$
and $H$ are indecomposable, then we apply
Lemma~\ref{lem:trivialJSJ}. If both $G$ and $H$ have non-trivial
Grushko decompositions with the same number of factors, then by
Lemma~\ref{lem:free decomposition}, it is enough to compare the
factors $G_i$ and $H_j$ of these decompositions. If $G_i$ and
$H_j$ are free groups, then they are isomorphic if and only if
their generating sets have the same cardinality. Otherwise, $G_i$
and $H_j$ are one-ended groups (in what follows, we still denote
these groups by $G$ and $H$), and we consider their Abelian JSJ
decompositions $\Gamma(V,E)$ and $\Xi(U,P)$.
Theorem~\ref{thm:down_to_vertices} gives rise to the following
algorithm. We find all possible isomorphisms between the graphs
$(V,E)$ and $(U,P)$. If there are not any, then we are done as the
groups are not isomorphic. Otherwise, fix an isomorphism
$\lambda\colon (V,E)\rightarrow(U,P)$ and try to find an
extendable isomorphism $\varphi_v\colon G_v\rightarrow
H_{\lambda(v)}$ that preserves the ordering of the edge subgroups
(see the beginning of this section), for each $v\in V$. This
latter procedure depends on the type of the vertex group in
question: Abelian (elementary), flexible or rigid non-elementary.
Recall that by Lemma~\ref{lem:classify}, we are able to determine
the type of each vertex group effectively. If $G_v$ and
$H_{\lambda(v)}$ are either elementary or flexible groups, then it
suffices to compare their canonical presentations that are output
of the Elimination process. The groups are isomorphic if and only
if a map sending the generators of $G_v$ in the canonical
presentation to the generators of $H_{\lambda(v)}$, sends the
peripheral subgroups of $G_v$ onto the peripheral subgroups of
$H_{\lambda(v)}$, so that it remains to check that the ordering of
the peripheral subgroups is preserved. An algorithm for rigid
groups is the content of Lemma~\ref{lem:rigid gps}. Observe that
each rigid non-elementary subgroup is an $\mathcal{F}$-group with
the trivial Abelian decomposition modulo the set of peripheral
subgroups, which makes Lemma~\ref{lem:rigid gps} applicable in
this case. If for each isomorphism of graphs $\lambda$ there is a
pair of vertices $(v,\lambda(v))$ with no e-isomorphism between
$G_v$ and $H_{\lambda(v)}$ preserving the ordering (which we can
find out in a finite time), then $G$ and $H$ are not isomorphic.
Otherwise, we fix $\lambda$ and an e-isomorphism $\varphi_v$ for
each pair $(v,\lambda(v))$ and associate the comparative labelling
as defined above, to each set of e-isomorphisms between the
non-elementary vertices of $(V,E)$ and $(U,P)$. Since by
Corollary~\ref{cor:BPproperty}, the set of e-isomorphisms between
two rigid vertices is finite, we can apply
Lemma~\ref{lem:criterion} and obtain the claim.
\end{proof}
%
\section{Structure of the automorphism group}\label{s:autos}
Let $G$ be a one-ended $\mathcal{F}$-group. By
Theorem~\ref{thm:isomorphic image of JSJ}, an Abelian JSJ
decomposition of $G$ and its image under an automorphism of $G$
differ by conjugation and modifying boundary monomorphisms. We
apply this result to study the structure of $Out(G)$. To state our
result, we introduce one more definition.

\begin{defn}\label{defn:canonical autos} Let $G$ be a freely
indecomposable $\mathcal{F}$-group, and let $\Gamma(V,E\cup
E_s;T)$ be the Abelian JSJ decomposition of $G$. We define the
group $Out_{\Gamma}(G)$ to be the subgroup of $Out(G)$ generated
by the following types of automorphisms of $G$:
\begin{enumerate}
  \item \label{e:Dehn twist} Generalized Dehn twists along edges in $\Gamma$
  (see Definition~\ref{defn:Dehn twist}).
    \item \label{e:abelian vertices} Automorphisms of an elementary
  vertex group that preserve the peripheral subgroups of the
  group.
  \item \label{e:mapping class groups} Automorphisms of a flexible
  vertex group $G_u$ that preserve the peripheral subgroups of the
  group, up to conjugacy (geometrically, these are Dehn twists
  along simple closed curves on the punctured surface
  $\Sigma$ with $\pi_1(\Sigma)\cong G_u$).
  \end{enumerate}
\end{defn}
\begin{lem} With the notation of Definition~\ref{defn:canonical
autos}, $[Out(G):Out_{\Gamma}(G)]<\infty$.
\end{lem}
\begin{proof} According to Theorem~\ref{thm:isomorphic image of JSJ}, each automorphism
$\psi\in Aut(G)$ preserves the maximal tree $T$ of $\Gamma$.
Therefore, $\psi$ is the composition of e-automorphisms of
vertices, automorphisms of type~\refe{Dehn twist}, and
conjugation. Observe that the e-automorphisms of elementary and
flexible vertices belong to $Out_{\Gamma}(G)$. Furthermore,
according to Corollary~\ref{cor:BPproperty}, each rigid vertex has
only finitely many e-automorphisms. Also observe that
e-automorphisms of different vertices commute; the assertion
follows.
\end{proof}

Let $G$ be a one-ended $\mathcal{F}$-group, and let $\Gamma(V,E)$
be an Abelian JSJ decomposition of $G$. By an \emph{e-automorphism
of a vertex group} $G_v$ we mean an automorphism $\psi\in
Out(G_v)$ that maps each edge subgroup of $G_v$ onto a conjugate
of itself (cf. Definition~\ref{defn:extendable iso}). We denote by
$V_M\subset V$ the subset of all elementary vertices and by
$V_Q\subset V$ the subset of all flexible (or QH-)vertices of
$\Gamma$ (see~\cite[Definition]{RS2} and
Definition~\ref{defn:QHsbgp} in the present paper). With each
vertex $v\in V_M\cup V_Q$ we associate the subgroup of
e-automorphisms of $G_v$ denoted by $\mathcal{M}_v$ if $v\in V_M$
and by $\mathcal{Q}_v$ if $v\in V_Q$. Since $G_v$ is a finitely
generated free Abelian group, $\mathcal{M}_v$ is a subgroup of
$GL_n(\mathbb{Z})$, where $n$ is the maximal rank of an Abelian
subgroup of $G$. Each flexible vertex group is the fundamental
group of a punctured surface, so that $\mathcal{Q}_v$ is the
mapping class group of a surface with boundary. Let
$\mathcal{Q}=\prod_{v\in V_Q}\mathcal{Q}_v$ and
$\mathcal{M}=\prod_{v\in V_M}\mathcal{M}_v$. Since the structure
of $Out_{\Gamma}(G)$ is well understood, we have the following
result.
\begin{thm}\label{thm:auto group} Let $G$ be a one-ended
$\mathcal{F}$-group. The group $Out(G)$ is virtually a direct
product $\mathbb{Z}^d\times\mathcal{M}\times\hat{\mathcal{Q}}$
where $\hat{\mathcal{Q}}$ is the quotient of $\mathcal{Q}$ by a
central subgroup isomorphic to a f.g. free Abelian group
$\mathbb{Z}^m$.
\end{thm}

\end{document}